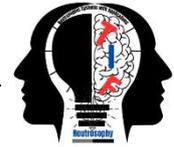

# Pura Vida Neutrosophic Algebra

**Ranulfo Paiva Barbosa (Sobrinho)** [1,*] and **Florentin Smarandache** [2]

[1] Web3 Blockchain Entrepreneur; 37 Dent Flats, Monte de Oca, 11501, San José, Costa Rica; ranulfo17@gmail.com.
[2] Math and Science Department, University of New Mexico, Gallup, NM, 87301, USA; smarand@unm.edu.

* Correspondence: ranulfo17@gmail.com.

**Abstract:** We introduce Pura Vida Neutrosophic Algebra, an algebraic structure consisting of neutrosophic numbers equipped with two binary operations namely addition and multiplication. The addition can be calculated sometimes with the function min and other times with the max function. The multiplication operation is the usual sum between numbers. Pura Vida Neutrosophic Algebra is an extension of both Tropical Algebra (also known as Min-Plus, or Min-Algebra) and Max-Plus Algebra (also known as Max-algebra). Tropical and Max-Plus algebras are algebraic structures included in semirings and their operations can be used in matrices and vectors. Pura Vida Neutrosophic Algebra is included in Neutrosophic semirings and can be used in Neutrosophic matrices and vectors.

## 1. Introduction

Uncertain, indeterminacy, imprecise, and vague are common characteristics of data in real-life problems like decision-making, engineering, computer science, finance, etc. Several theories have been proposed to deal with these data characteristics, fuzzy set theory [1], intuitionistic fuzzy sets [2], rough set theory [3], Soft set [4], and Neutrosophy theory [5]. Since Smarandache introduced Neutrosophy to study the basis, nature, and range of neutralities as well as their contact with ideational spectra in the 1990s, we have seen the emergence of neutrosophic algebraic structures [6], neutrosophic probability and statistics [7, 8] neutrosophic numbers [8], single-valued neutrosophic sets (SVNSs) [9, 21], and several algebraic structures such as neutrosophic semirings [10], among others theoretical advances [11] and also applications [12].

Through neutrosophic semirings, we introduce Pura Vida (PV) Neutrosophic Algebra, an algebraic structure consisting of neutrosophic numbers equipped with two binary operations namely addition and multiplication. Pura Vida Neutrosophic Algebra is an extension of both Tropical Algebra (also known as Min-Plus) [13] and Max-Plus Algebra [14]. Both Tropical and Max-Plus algebra are algebraic structures included in semirings and were discovered independently by several researchers [13, 14]. They were defined on the real number domain and for the first time, we extended them to the neutrosophic domain.

## 2. Preliminaries
### 2.1. Semiring

A semiring [15] denoted (V, $\oplus$, $\otimes$, 0, 1 ) is a set V equipped with two binary operations, addition:

$$\oplus : V \times V \rightarrow V$$

And multiplication:

$$\otimes: V \times V \to V$$

Which satisfies the following axioms for any u, v, w ∈ V:

1. (V, ⊕, 0) is a commutative monoid and (V, ⊗, 1) is a monoid.
2. u ⊗ (v ⊕ w) = (u ⊗ v) ⊕ (u ⊗ w) and (v ⊕ w) ⊗ u = (v ⊗ u) ⊕ (w ⊗ u) (distributivity).
3. 0 annihilates V: v ⊗ 0 = 0 ⊗ v = 0.

When (V, ⊗, 1) is a commutative monoid, the semiring (V, ⊕, ⊗, 0, 1 ) is said to be a commutative semiring.

## 2.2 Tropical algebra

The tropical algebra is also referred to as tropical semiring T, which consists of the set of real numbers, R, extended with infinity, equipped with the operations of taking minimums (as semiring addition) and addition (as semiring multiplication) [14, 16]. Tropical algebra is also known as min-plus algebra. With minimum replaced by maximum, we get the isomorphic max-plus algebra [17]. According to [17], the adjective "tropical" was coined by French mathematicians to honor their Brazilian colleague Imre Simon [16], who pioneered the use of min-plus algebra in optimization theory.

$$T = (R \cup \{ \infty \}, \oplus, \otimes )$$

Addition operation:
$$a \oplus b = \min(a, b)$$

Multiplication operation:
$$a \otimes b = a + b$$

the operations of R, are extended to T in the usual way and the identities of ⊕ and ⊗ are, respectively, ∞ and 0. The element ∞ represents plus-infinity [13]. Given a real number, x ∈ T, its addition and multiplication identity are given, respectively:

$$x \oplus \infty = x$$
$$x \otimes 0 = x$$

Michaleck points out the following equations involving the two identity elements:

$$x \otimes \infty = \infty \quad \text{and} \quad x \oplus 0 = \begin{cases} 0, & \text{if } x \geq 0 \\ x, & \text{if } x < 0 \end{cases}$$

Michaleck said there is no subtraction in tropical arithmetic. Tropical division ⊘ is defined to be classical subtraction.

Tropical division, x ⊘ y = x, exists if and only if y ⊗ z = x [20].

In Tropical algebra the pairs of operations (⊕,⊗) is extended to matrices and vectors similarly as in linear algebra. That is if A = ($a_{ij}$), B = ($b_{ij}$) and C = ($c_{ij}$) are matrices with elements from R of compatible sizes, we write:

C = A ⊕ B if $c_{ij}$ = $a_{ij}$ ⊕ $b_{ij}$ for all i, j
C = A ⊗ B if $c_{ij}$ = $\sum_{k}^{\oplus} a_{ik} \oplus b_{kj}$ = $\max_k (a_{ik} + b_{kj})$ for all i, j
$\alpha$ ⊗ A = A ⊗ $\alpha$ = ( $\alpha$ ⊗ $a_{ij}$ ) for all $\alpha$ ∈ R.

## 2.3 Max-Plus algebra

The Max-Plus algebra is an algebraic structure semiring MP, which consists of the set of real numbers, R, extended with infinity, equipped with the operations of taking maximums (as semiring addition) and addition (as semiring multiplication) [14].

$$MP = (R \cup \{ -\infty \}, \oplus', \otimes)$$

Addition operation:
$$a \oplus' b = \max(a, b)$$

Multiplication operation:
$$a \otimes b = a + b$$

the operations of R, are extended to MP in the usual way and the identities of $\oplus'$ and $\otimes$ are, respectively, $-\infty$ and 0.

In max-plus algebra the pairs of operations ($\oplus'$, $\otimes$) is extended to matrices and vectors similarly as in linear algebra. That is if A = ($a_{ij}$), B = ($b_{ij}$) and C = ($c_{ij}$) are matrices with elements from R of compatible sizes, we write:

C = A $\oplus'$ B if $c_{ij}$ = $a_{ij}$ $\oplus'$ $b_{ij}$ for all i, j

C = A $\otimes$ B if $c_{ij}$ = $\sum_{k}^{\oplus'} a_{ik} \oplus b_{kj}$ = $\max_{k}(a_{ik} + b_{kj})$ for all i, j

$\alpha \otimes A = A \otimes \alpha = (\alpha \otimes a_{ij})$ for all $\alpha \in$ R.

## 2.4 Neutrosophic Set

Smarandache [5] defined Neutrosophic set as a set of elements composed of tripart structure: a Truth membership (T), an Indeterminacy membership (I) and a False membership (F). These parts are independent each other and can be represented by different functions. Together, <T, I, F>, these parts compose an element of Neutrosophic set.

## 2.5 Neutrosophic Number

According to [18] the neutrosophic number (NN) is a number which structure is given by "X = a + bI", where I represents the indeterminacy component of X, and 'a' and 'b' are real or complex numbers [19].

## 2.6 Neutrosophic Semiring

An algebraic structure (S∪I, $\oplus$, $\otimes$) is called neutrosophic semiring [10] if $\oplus$ and $\otimes$ are the closed and associative binary operations and $\otimes$ is distributive over $\oplus$, where S is semiring with respect to $\oplus$ and $\otimes$ and I is the neutrosophic element (I = $I^2$) and < S∪I> = { a + bI; a, b $\in$ S}.

## 2.7 Neutrosophic field [6]

Let K be the field of reals. We call the field generated by K $\cup$ I to be the neutrosophic field for it involves the indeterminacy (I) factor in it. We define $I^2$ = I, I + I = 2I, i.e., I + … + I = nI, and if k $\in$ K then kI = Ik, 0I = 0. We denote the neutrosophic field by K(I).

## 2.8 Neutrosophic matrix [6]

Let $M_{nxm}$ = {($a_{ij}$) / $a_{ij}$ $\in$ K(I) }, where K(I), is a neutrosophic field. We call Mnxm to be the neutrosophic matrix.

## 3. Pura Vida Neutrosophic Algebra

The Pura Vida Neutrosophic Algebra, PV, is an extension of the Tropical algebra and Max-Plus Algebra.

Pura Vida Neutrosophic Algebra is included in a Neutrosophic semiring, i.e., it has both associative binary operations, addition $\oplus$ and multiplication $\otimes$ where $\otimes$ is distributive over $\oplus$, and S is semiring with respect to $\oplus$ and $\otimes$ and I is the neutrosophic element (I = $I^2$) and < S∪I> = { a = bI; a, b $\in$ S}. The addition operation can use either the min function, $\oplus$, or the max function, $\oplus'$, depending on the situation.

$$PV = (S \cup I \{-\infty, +\infty\}, \oplus, \oplus', \otimes)$$

Pura Vida Neutrosophic Algebra operations addition ($\oplus$, or, $\oplus'$) and multiplication ($\otimes$) are given:

## 3.1 Addition operation $\oplus$, or, $\oplus'$

Depending on the real-life applications, the addition operation can use the min or max function. Given two neutrosophic numbers x = a + bI, and z = c + dI $\in$ S, the addition of x and z:

**3.1.1.** $\quad$ x $\oplus$ z = ( a $\oplus$ c ) + ( b $\oplus$ d )I = min(a, c) + min(b, d)I

or,

**3.1.2.**    $x \oplus' z = (a \oplus' c) + (b \oplus' d)I = \max(a, c) + \max(b, d)I$

## 3.2 Multiplication operation $\otimes$

Given two neutrosophic numbers $x = a + bI$, and $z = c + dI \in S$, the multiplication of x and z:

   $x \otimes z = a \otimes c + (b \otimes d)I = (a + c) + (b + d)I$

## 3.3 Identities

In Pura Vida Neutrosophic Algebra, PV, the identities of the operators $\oplus$, $\oplus'$ and $\otimes$ are, respectively, $\infty$, $-\infty$ and 0.

## 3.4 Properties

Next, we show that the PV attends the closure property and distributive and associative laws. We use min for the addition operation, but, one could use the max function to show that PV verifies the mentioned properties.

### 3.4.1 Closure property:

Let $(a + bI)$ and $(c + dI) \in S \cup I$ then,

   $(a + bI) \oplus (c + dI) = (a \oplus c) + (b \oplus d)I = \min(a, b) + \min(c, d)I, \in S \cup I$. The addition operation verifies the closure property.

   $(a + bI) \otimes (c + dI) = a \otimes c + (b \otimes d)I = (a + c) + (b + d)I \in S \cup I$. Which shows that the closure property is satisfied for the multiplication operation.

### 3.4.2 Distributive law:

Let $(a + bI)$, $(c + dI)$ and $(e + fI) \in S \cup I$, then:

   $(a + bI) \otimes [(c + dI) \oplus (e + fI)] = (a + bI) \otimes [\min(c, e) + \min(d, f)I] =$
   $= [a + \min(c, e)] + [b + \min(d, f)]I$.
   And $[(a + bI) \otimes (c + dI)] \oplus [(a + bI) \otimes (e + f)I] =$
   $= [(a + c) + (b + d)I] \oplus [(a + e) + (b + f)I] =$
   $= \min\{(a + c), (a + e)\} + \min\{(b + d) + (b + f)\}I =$
   $= [a + \min(c, e)] + [b + \min(d, f)]I$.

### 3.4.3 Associative law:

Let $(a + bI)$, $(c + dI)$ and $(e + fI) \in S \cup I$, then:

   $[(a + bI) \oplus (c + dI)] \oplus (e + fI) =$
   $[\min(a, c) + \min(b, d)I] \oplus (e + fI) = \min[\min(a, c), e] + \min[\min(b, d), f]I =$
   $= (a \oplus c \oplus e) + (b \oplus d \oplus f)I$.
   $(a + bI) \oplus [(c + dI) \oplus (e + fI)] =$
   $(a + bI) \oplus [\min(c, e) + \min(d, f)I] = \min[a, \min(c, e)] + \min[b, \min(d, f)]I =$
   $= (a \oplus c \oplus e) + (b \oplus d \oplus f)I$.
   Again:
   $[(a + bI) \otimes (c + dI)] \otimes (e + fI) = [a \otimes c + (b \otimes d)I] \otimes (e + fI) =$
   $[(a + c) + (b + d)I] \otimes (e + fI) = (a + c) \otimes e + [(b + d) \otimes f]I =$
   $(a + c + e) + (b + d + f)I$.
   $(a + bI) \otimes [(c + dI) \otimes (e + fI)] = (a + bI) \otimes [(c + e) + (d + f)I] =$
   $= a \otimes (c + e) + b \otimes (d + f)I = (a + c + e) + (b + d + f)I$.

## 3.5 Pura Vida Neutrosophic Algebra on Matrices

In Pura Vida Neutrosophic Algebra the pairs of operations $(\oplus, \oplus', \otimes)$ is extended to matrices and vectors similarly as in linear algebra. That is if $A = (a_{ij})$, $B = (b_{ij})$ and $C = (c_{ij})$ are matrices with elements from R of compatible sizes, we write:

   $C = A \oplus' B$ if $c_{ij} = a_{ij} \oplus' b_{ij}$ for all i, j
   $C = A \otimes B$ if $c_{ij} = \sum_k^{\oplus'} a_{ik} \oplus b_{kj} = \max_k (a_{ik} + b_{kj})$ for all i, j

$\alpha \otimes A = A \otimes \alpha = (\alpha \otimes a_{ij})$ for all $\alpha \in R$.

### 3.5.1 Matrices Addition using $\oplus$ operator

Given P and Q, both square neutrosophic matrices 2x2, their sum is $D = P \oplus Q$.

$$P = \begin{vmatrix} -8+I & 5-I \\ 3+8I & 23-2I \end{vmatrix} \text{ and } Q = \begin{vmatrix} 3+2I & 13+3I \\ 7+9I & 3+5I \end{vmatrix}$$

$$D = \begin{vmatrix} \text{Min}(-8,3)+\text{Min}(1,2)I = -8+I & \text{Min}(5,13)+\text{Min}(-1,3)I = 5-I \\ \text{Min}(3,7)+\text{Min}(8,9)I = 3+8I & \text{Min}(23,3)+\text{Min}(-2,5)I = -2+3I \end{vmatrix}$$

### 3.5.2 Matrices Addition using $\oplus'$ operator

Given X and Z, both square neutrosophic matrices 2x2, their sum is $W = X \oplus' Z$.

$$X = \begin{vmatrix} -8+I & 5-I \\ 3+8I & 23-2I \end{vmatrix} \text{ and } Z = \begin{vmatrix} 3+2I & 13+3I \\ 7+9I & 3+5I \end{vmatrix}$$

$$W = \begin{vmatrix} \text{Max}(-8,3)+\text{Max}(1,2)I = 3+2I & \text{Max}(5,13)+\text{Max}(-1,3)I = 13+3I \\ \text{Max}(3,7)+\text{Max}(8,9)I = 7+9I & \text{Max}(23,3)+\text{Max}(-2,5)I = 23+5I \end{vmatrix}$$

### 3.5.3 Matrices Multiplication using $\otimes$ operator

Given A and B, both rectangular neutrosophic matrices, their multiplication is $C = A \otimes B$.

$$A = \begin{vmatrix} -1 & 2 & -I \\ 3 & I & 0 \end{vmatrix} \text{ and } B = \begin{vmatrix} I & 1 & 2 & 4 \\ 1 & I & 0 & 2 \\ 5 & -2 & 3I & -I \end{vmatrix}$$

$$C = A \otimes B = \begin{vmatrix} C_{11} & C_{12} & C_{13} & C_{14} \\ C_{21} & C_{22} & C_{23} & C_{24} \end{vmatrix}$$

Where,

$C_{11} = (-1 \ 2 \ -I) \otimes (I \ 1 \ 5) = (-1 \otimes I + 2 \otimes 1 + -I \otimes 5) = (-1+I + 2+1 + -I+5) = 7$

$C_{21} = (3 \ I \ 0) \otimes (I \ 1 \ 5) = (3 \otimes I + I \otimes 1 + 0 \otimes 5) = (3+I + 1+I + 5) = 9+2I$.

$C_{12} = (-1 \ 2 \ -I) \otimes (1 \ I \ -2) = (-1 \otimes 1 + 2 \otimes I + -I \otimes -2) = (0 + 2 + I -I -2) = 0$.

$C_{22} = (3 \ I \ 0) \otimes (1 \ I \ -2) = 3 + 1 + I + I -2 = 2 + 2I$.

$C_{13} = (-1 \ 2 \ -I) \otimes (2 \ 0 \ 3I) = 1 + 2 + 2I = 3 + 2I$.

$C_{23} = (3 \ I \ 0) \otimes (2 \ 0 \ 3I) = 5 + I + 3I = 5 + 4I$.

$C_{14} = (-1 \ 2 \ -I) \otimes (4 \ 2 \ -I) = 3 + 4 -2I = 7 - 2I$.

$C_{24} = (3 \ I \ 0) \otimes (4 \ 2 \ -I) = 7 + I + 2 -I = 9$.

$$C = A \otimes B = \begin{vmatrix} 7 & 0 & 3+2I & 7-2I \\ 9+2I & 2+2I & 5+4I & 9 \end{vmatrix}$$

## 4. Conclusion

We introduced Pura Vida Neutrosophic Algebra through neutrosophic numbers and explored some its properties and applied to neutrosophic matrices.

**Data availability**

The datasets generated during and/or analyzed during the current study are not publicly available due to the privacy-preserving nature of the data but are available from the corresponding author upon reasonable request.

**Conflict of interest**

The authors declare that there is no conflict of interest in the research.

**Ethical approval**

This article does not contain any studies with human participants or animals performed by any of the authors.